\documentclass[11pt]{amsart}%
\usepackage{amsmath}
\usepackage{graphicx}
\usepackage{amsfonts}
\usepackage{amssymb}
\usepackage{hyperref}%
\usepackage{xcolor}

\usepackage{colortbl}

\providecommand{\U}[1]{\protect\rule{.1in}{.1in}}

\theoremstyle{plain}

\newtheorem{definition}{Definition}[section]

\newtheorem{lemma}{Lemma}[section]

\newtheorem{remark}{Remark}[section]

\numberwithin{equation}{section}

\title[Intrinsicality of second fundamentals form]{Intrinsicality of second fundamentals of hypersurfaces in space form\\ }
\dedicatory{Dedicated to the memory of Professor Joseph Kohn}
\author{Pengfei Guan}
\address{Department of Mathematics and Statistic\\
McGill University, Montreal, Quebec,\ Canada}
\email{pengfei.guan@mcgill.ca}
\author{Yu Yuan}
\address{Department of Mathematics, Box 354350\\
University of Washington\\
Seattle, WA 98195}
\email{yuan@math.washington.edu}
\date{\today}

\begin{document}

\begin{abstract}
 We prove that the normalized second fundamental form $A$ of immersed $C^2$ hypersurface in a space form $(N^{n+1}, \bar g)$ is intrinsic provided $\sigma_{2k+1}(A)\neq 0$ for some $k\ge 1$. We also establish the intrinsicality of the normalized second fundamental form $A$ of $M^{n}=\partial\Omega$ for domain $ \Omega\subset N^{n+1}, \ n\ge 3$.
\end{abstract}

\maketitle

\section{Introduction}

Gauss' Theorema Egregium \cite{Gauss} states that Gauss curvature of any immersed surface $M$ in $\mathbb R^3$ is intrinsic. For higher dimensions, Weyl \cite{Weyl} established that $\sigma_{2k}(A)$ is intrinsic for immersed hypersurface $M\in \mathbb R^{n+1}$, where $A=g^{-1}h$, $g$ is the induced metric, $h$ is the 
second fundamental form of $M,$ and $\sigma_m$ is the $m$-th elementary symmetric function. 

In this note, we discuss the intrinsicality of the second fundamental form of hypersurfaces in space form $N^{n+1},\ n\ge 3$.  We establish the intrinsicality of odd symmetric combinations of
principal curvatures of hypersurfaces in space forms under certain rank conditions. This rank condition is not needed when $M$ is the boundary of a bounded domain in $N^{n+1}, \ n\ge 3$. The proof is simple and elementary.

\medskip

From now on,  $\left(  \mathbb{N}^{n+1}, \bar g\right)$ denotes a space form with constant curvature $\bar{K}=0,\pm1$. 

\begin{definition} A geometric quantity $Q$ on a Riemannian manifold $(M,g)$ is called intrinsic if there is a continuous function $\phi$ such that
\[Q(x)=\phi(Riem_g(x)), \quad \forall x\in M,\] where $Riem_g=R^{\alpha\beta}_{\gamma\eta}$ is the {\color{blue} $(2,2)$} Riemann tensor of $(M,g).$
\end{definition}

The first result is of local nature.
\begin{lemma}\label{l1}
Let $A=g^{-1}h$ be the normalized second fundamental form of an immersed $C^2$ hypersurface $\left(
M^{n},g\right)  \subset\left(  \mathbb{N}^{n+1}, \bar g\right)  $ with the intrinsic induced metric $g =\left.  \bar g\right\vert _{M^{n}},$ then
\begin{enumerate}\item $\sigma_{2k+1}\left(  A\right)  \sigma_{2l+1}\left(  A\right), \ \ \forall\  k+l\geq1$is intrinsic. More precisely, for given pair $(2k+1,$ $2l+1)$, there is a polynomial $P_{2k+1,2l+1}$ of degree $k+l+1$ in term of $(2,2)$ Riemnnian tensor $Riem_g$ such that
\begin{equation}\label{intrin} \sigma_{2k+1}\left(  A\right)  \sigma_{2l+1}\left(  A\right)=P_{2k+1,2l+1}(Riem_g).\end{equation}
\item If $\sigma_{2k+1}\left(  A\right)  \neq0$ for some
$k\geq1$, then $\sigma_{2l+1}
$ is intrinsic $\forall\  l$ modulo orientation. As a consequence, $A$ is intrinsic up to orientation. \item If the rank of $A$ is at least three, $H^2=\sigma_{1}^{2}\left(  A\right)  $ is intrinsic.  If in addition $H\neq 0$, then  $A$ is intrinsic up to orientation.\end{enumerate}
\end{lemma}

\medskip
The next is of global nature.
\begin{lemma}\label{l3}
Let   $\Omega\subset\left(\mathbb N^{n+1},\bar g, \bar K=0\  \text{or}\ -1\right)$ or  $\left(\mathbb S^{n+1}_{+}, \bar g, \  \bar K=1\right)$,  $\ n\ge 3$ be a domain with $C^2$ connected and compact boundary $M^{n}=\partial\Omega$ with outward orientation,  then the normalized second fundamental form $A$ is intrinsic.   

As a consequence,
\[\int_{M}\sigma_{k}^{m}\left(  A\right)  d\sigma\] is intrinsic for all
$k=0,1,\cdots,n$ and all positive integers $m.$ 

\end{lemma}


Lemma \ref{l3} extends Weyl's even results \cite{Weyl} in the case of compact hypersurfaces.

\section{Proof }
\subsection{Proof of Lemma \ref{l1}}
\begin{proof}
At a fixed point, we assume $A$ is diagonal with an orthonormal frame $e_1, e_2, \cdots, e_n$ on $M$. The eigenvalues of $A$ are the principal
curvatures of the hypersurface $\kappa_{1},\cdots,\kappa_{n}.$ 
By Gauss equation, $\forall\ \alpha\neq \beta$,
\begin{equation}\label{Gauss} \kappa_{\alpha}\kappa_{\beta}=R_{\alpha\beta\alpha\beta}-\bar{K},\end{equation} where $R_{\alpha\beta\alpha\beta}$ is the sectional curvature with respect to the frame $\{e_1, e_2, \cdots, e_n\}$ of the intrinsic metric $g$.

Observe that the $C^{2}$ smoothness of the immersion $\left(  M^{n},g\right)
$ in the space form is enough to determine its Riemann curvature, which is in
terms of the derivatives of the metric $g$ up to second order, though the
induced metric appears only $C^{1}.\ $In fact, near each point $o$ of the
immersion, we parametrize by its tangent coordinates $\left(  x,u\left(
x\right)  \right)  $ with $Du\left(  o\right)  =0.$ This can be realized by a
(Euclid) rotation in the Riemann model of the space form metric, coupled with
the implicit function theorem. The induced metric takes the following form%
\begin{gather*}
\left(  g_{ij}\right)  =\frac{4}{  1-\bar{K}\left\vert X\right\vert
^{2}  }\left.  dX^{2}\right\vert _{M}=\left(  \frac{4\left[
\delta_{ij}+u_{i}\left(  x\right)  u_{j}\left(  x\right)  \right]
dx_{i}dx_{j}}{1-\bar{K}\left\vert \left(  x,u\left(  x\right)  \right)
\right\vert ^{2}}\right)  .
\end{gather*}
Because of $Du\left(  o\right)  =0,$ one sees that the second order derivatives of the $C^{2}$ intrinsic metric $g$ only depend on
up to second order derivatives of the immersion for the Riemann curvature at
$o$%
\[
R_{ijkl}\left(  o\right)  =\frac{1}{2}\left(  g_{il,jk}+g_{jk,il}%
-g_{jl,ik}-g_{ik,jl}\right).\]

From the above discussion, if $M=\{F=0\}$  locally with $\nabla F\neq 0$, then the Riemann tensor $Riem_g$ is a continuous function of $(\nabla F, \nabla^2F)$.
\medskip

\begin{enumerate}
\item We assume $k\geq l,$ then $k\geq1.$ 
We write
\[
\sigma_{2k+1}\left(  A\right)  \sigma_{2l+1}\left(  A\right)  =\sigma
_{2k+1}\left(  \kappa\right)  \sigma_{2l+1}\left(  \kappa\right)  =\frac
{1}{2k+1}\sum_{i=1}^{n}\sigma_{2k}\left(  \left.  \kappa\right\vert i\right)
\kappa_{i}\sigma_{2l+1}\left(  \kappa\right)  .
\]
For each $i,$%
\[
\sigma_{2k}\left(  \left.  \kappa\right\vert i\right)  \kappa_{i}\sigma
_{2l+1}\left(  \kappa\right)  =\sigma_{2k}\left(  \left.  \kappa\right\vert
i\right)  \left[  \sigma_{2l+1}\left(  \left.  \kappa\right\vert i\right)
\kappa_{i}+\sigma_{2l}\left(  \left.  \kappa\right\vert i\right)  \kappa
_{i}^{2}\right]  .
\]
Note that $\sigma_{2k}\left(  \left.  \kappa\right\vert i\right)  $ is a
combination of $\kappa_{\alpha}\kappa_{\beta}$ for $\alpha\neq\beta\neq i,$
and also $\sigma_{2l+1}\left(  \left.  \kappa\right\vert i\right)  \kappa_{i}$
is a combination of $\kappa_{\alpha}\kappa_{\beta}$ and $\kappa_{\gamma}%
\kappa_{i}$ for $\alpha\neq\beta\neq i$ and $\gamma\neq i.$ Then both
$\sigma_{2k}\left(  \left.  \kappa\right\vert i\right)  \ $and $\sigma
_{2l+1}\left(  \left.  \kappa\right\vert i\right)  \kappa_{i}$ can be
expressed as a function of section curvatures $R_{\alpha\beta\alpha\beta}$ and
$R_{\gamma i\gamma i}$ of the metric $g.$ Moreover, for $k\geq1,$ $\sigma
_{2k}\left(  \left.  \kappa\right\vert i\right)  \kappa_{i}^{2}$ can be
expressed in terms of $\kappa_{\alpha}\kappa_{\beta},\kappa_{\alpha}\kappa
_{i},$ and $\kappa_{\beta}\kappa_{i},$ which are all sectional curvatures of
metric $g,$  for $\alpha\neq\beta,$ $\alpha\neq i,$ and $\beta\neq i.$

Summing up, one sees that $\sigma_{2k+1}\left(  \kappa\right)  \sigma
_{2l+1}\left(  \kappa\right)  $ can be express as a polynomial $P_{P_{2k+1,2l+1}}(Riem_g)$ of degree $k+l+1$ in term of Riemannian tensor $Riem_g$. In particular, $\sigma_{2k+1}^2(A)$ is intrinsic,
\[\sigma_{2k+1}^2(A)=P_{2k+1, 2k+1}(Riem_g).\]

\medskip
We remark that since $\sigma_{2k}(\kappa)  $ is a
combination of $\kappa_{\alpha}\kappa_{\beta}$ for $\alpha\neq\beta$. By Gauss equation (\ref{Gauss}) and the above argument, Weyl's intrinsicality of $\sigma_{2k}$ holds true for general hypersurfaces in space forms. 

\medskip

\item Consequently, if $\sigma_{2k+1}\left(  A\right)  \neq0$, 
by (\ref{intrin}),
\[\sigma_{2k+1}\left(  A\right) =\pm{|P_{2k+1,2k+1}(Riem_g)|^{\frac12}},\] where $\pm$ sign is determined by orientation. Once the sign for $\sigma_{2k+1}(A)$ is given, $\sigma_{2l+1}(A), \ \forall \  l$ is completely determined by (\ref{intrin}) and so is the orientation. Of course, the reverse is also true.  Moreover
\[
\sigma_{2l+1}\left(  A\right) 
=\frac{\left[\sigma_{2l+1}\left(  A\right)  \sigma_{2k+1}\left(  A\right) \right] }%
{\sigma_{2k+1}\left(  A\right)  }=\pm{\frac{P_{{2k+1,2l+1}}(Riem_g)}{|P_{2k+1,2k+1}(Riem_g)|^{\frac12}}}.%
\]
is also intrinsic up to orientation. 

Combined with the well-known fact that an even $\sigma_{2k}\left(  A\right)  $
is intrinsic unconditionally, which can be shown by a simple even arrangement
of $\sigma_{2k}\left(  A\right)  ,$ one sees all principal curvatures
$\kappa_{1},\cdots,\kappa_{n}$ are intrinsic up to orientation, and so is $A.$

\medskip

\item We prove the last statement of the Lemma. One has \[\sigma_{1}^{2}\left(  \kappa\right)  =\kappa_{1}^{2}+\cdots
+\kappa_{n}^{2}+2\sigma_{2}\left(  \kappa\right)  ,\] as the scalar curvature
$\sigma_{2}\left(  \kappa\right)  $ is intrinsic, we only need to show
$\left\vert A\right\vert ^{2}=\left\vert \kappa\right\vert ^{2}$  is intrinsic
under assumption that the rank of $A\ $is at least three. 

\medskip

If the rank of $A$ is odd $r=2k+1\geq3,$ then $\sigma_{r}\left(
\kappa\right)  \neq0$ and%
\[
\left\vert \kappa\right\vert ^{2}=\frac{1}{\sigma_{r}^{2}\left(
\kappa\right)  }\sum_{i=1}^{n}\left[  \kappa_{i}\sigma_{r}\left(
\kappa\right)  \right]  ^{2}.
\]
Observe that, as $r\geq3,$  $\kappa_{i}\sigma_{r}\left(  \kappa\right)  $ can
be expressed in terms of product of distinct pairs $\kappa_{i}\kappa_{\alpha
},\ \kappa_{i}\kappa_{\beta},$ and $\kappa_{\gamma}\kappa_{\delta},$ which are
all sectional curvature of the metric $g.\ $Combined with the intrinsicality
of $\sigma_{r}^{2}\left(  \kappa\right)  ,$ which is just proved in the above,
one sees the same intrinsicality for $\left\vert \kappa\right\vert ^{2}.$ 

\medskip

If the rank of $A$ is even $r=2l,$  then $l\geq2$ and $\sigma_{2l}\left(
\kappa\right)  \neq0.$  Consider%
\[
\sigma_{2l}\left(  \kappa\right)  \left\vert \kappa\right\vert ^{2}%
=\sum_{\alpha_{1}<\cdots<\alpha_{2l}}\sum_{j=1}^{n}\kappa_{j}^{2}%
\kappa_{\alpha_{1}}\cdots\kappa_{\alpha_{2l}}.
\]
In $\kappa_{j}^{2}\kappa_{\alpha_{1}}\cdots\kappa_{\alpha_{2l}},$ if
$\alpha_{i}\neq j,$ for all $i=1,\cdots,2l,$ then%
\[
\kappa_{j}^{2}\kappa_{\alpha_{1}}\cdots\kappa_{\alpha_{2l}}=\left(  \kappa
_{j}\kappa_{\alpha_{1}}\right)  \left(  \kappa_{j}\kappa_{\alpha_{2}}\right)
\left(  \kappa_{\alpha_{3}}\kappa_{\alpha_{4}}\right)  \cdots\left(
\kappa_{\alpha_{2l-1}}\kappa_{\alpha_{2l}}\right)
\]
is intrinsic, as $\kappa_{\alpha}\kappa_{\beta}$ is intrinsic for $\alpha
\neq\beta.$

If there is $\alpha_{i}=j,$ then $\alpha_{m}\neq j$ for $m\neq i.$ Since
$l\geq2,$ there are $2l-1\geq3$ many $\alpha_{m}$'s not being $j.$ Take three
from $\alpha_{m}$'s, say $\beta_{1},\beta_{2},\beta_{3},$ now%
\[
\kappa_{j}^{2}\kappa_{\alpha_{1}}\cdots\kappa_{\alpha_{2l}}=\left(  \kappa
_{j}\kappa_{\beta_{1}}\right)  \left(  \kappa_{j}\kappa_{\beta_{2}}\right)
\left(  \kappa_{j}\kappa_{\beta_{3}}\right)  \prod\limits_{\alpha_{m}%
\neq\alpha_{i},\beta_{1},\beta_{2},\beta_{3}}\kappa_{\alpha_{m}},
\]
where the last product $\prod$ term consists of $2l-4$ even many distinct
$\kappa_{\alpha_{m}}$'s. It follows that $\kappa_{j}^{2}\kappa_{\alpha_{1}%
}\cdots\kappa_{\alpha_{2l}}$ is intrinsic. In turn, as  $\sigma_{2l}\left(
\kappa\right)  \neq0,\ $ $\left\vert \kappa\right\vert ^{2}=\sigma_{2l}\left(
\kappa\right)  \left\vert \kappa\right\vert ^{2}/\sigma_{2l}\left(
\kappa\right)  $ is intrinsic. 

\medskip

Finally if the rank of $A$ is at least $3$ and $H\neq 0$, then $H$ is intrinsic and the orientation is determined by the sign of $H$. By (\ref{intrin}),
\[\sigma_{2k+1}(A)=\frac{P_{2k+1, 1}(Riem_g)}{H}, \ \forall \ k\ge 1 \] is intrinsic. Hence $A$ is intrinsic up to orientation.\end{enumerate} 

\medskip

Note that $P_{2k+1, 2l+1}(Riem_g)$ in (\ref{intrin}) is an invariant quantity independent of the choice of orthonormal frames 
and it is independent of orientation.  
\end{proof}

\subsection{Proof of Lemma \ref{l3}}
\begin{proof}

We first assume $M\ $is analytic. Since $M$ is compact and $n\ge 3$,  \[\Gamma=\left\{  \sigma_{3}\left(
A\right)  =0\right\}  \] is of measure $0$.
Otherwise, $\Gamma$ contains an n-dimensional open set in $M.$ By analyticity
and connectedness of $M,$ we know $\sigma_{3}\left(  A\right)  \equiv0$ on
$M.$ But this contradicts the fact that $A>0$ w.r.t. the inward normal at the
touching points of convex geodesic n-sphere envelope of $M$ in the Euclid,
hyperbolic, or half elliptic $\mathbb S^{n+1}_{+}$ ambient space form.

Denote \[\tilde M=M\smallsetminus\Gamma.\] $\tilde M$ is dense in $M$ by analyticity. 
$A$ is intrinsic in $\tilde M$ by Lemma \ref{l1} if we fix outward orientation.  For each $k\ge 0$ and for each point $x_0\in \Gamma$, as $M\in C^2$, $A\in C^0$ on $M$,
\[\sigma_{2k+1}(x_0)=\lim_{x\in \tilde M, x\to x_0} \sigma_{2k+1}(x). \] This implies that odd $\sigma_{2k+1}(A)$ is intrinsic. In turn, combined with the even intrinsicality of $\sigma_{2k}(A),$ $A$ is intrinsic when $M^n$ is analytic.

\medskip

For general case, we approximate $(M, g)$ by $(M_{\varepsilon}, g_{\epsilon})\in C^{\omega}$. We may write 
\[\Omega=\{ F<0\}, \quad \nabla F|_{M}\neq 0,\]
where $F$ is a defining function of $\Omega$. By Whitney's Approximation Theorem (Lemma 5 in \cite{Whitney}), $\exists \  F_{\epsilon}\in C^{\omega}$ such that 
\[\lim_{\epsilon\to 0} F_{\epsilon}=F, \quad \mbox{in} \ C^2,\]
uniformly in a closed set containing $\bar \Omega$. Denote 
\[\Omega_{\epsilon}=\{ F_{\epsilon}<0.\}\]
When $\epsilon$ is sufficient small, $\partial \Omega_{\epsilon}$ is $C^{\omega}$, and $\forall \  x\in \partial \Omega$ there is a unique $x_{\epsilon}\in \partial \Omega_{\epsilon}$ such that 
\[dist(y, x)\ge dist(x_{\epsilon}, x) \  \text{for all} \  y\in\partial \Omega_{\epsilon}\]
From the discussion in the proof of Lemma \ref{l1}, $Riem_g$ and $Riem_{g_{\epsilon}}$ are continuous functions of $(\nabla F, \nabla^2 F)$ and $(\nabla F_{\epsilon}, \nabla^2 F_{\epsilon} )$ respectively.
\[\lim_{\epsilon\to 0}\|Riem_{g_{\epsilon}}-Riem_{g}\|_{L^{\infty}}=0, \]
where $g_{\epsilon}$ the induced metric of embedded analytic hypersurface $M_{\epsilon}=\partial \Omega_{\epsilon}$. 

\medskip

Denote $A^{\epsilon}$ the normalized second fundamental form of $M_{\varepsilon}$. $A_{\epsilon}$ is intrinsic with respect to the outward orientation.  $\forall\ x\in M$, $x_{\epsilon}\in M_{\epsilon}$, 
\[A(x)=\lim_{\epsilon\to 0} A^{\epsilon}(x_{\epsilon})\] is intrinsic.
\end{proof}

\begin{remark} Lemma \ref{l1} provides an explicit intrinsic formula for $\sigma_{2l+1}(A), \ \forall \  l$ if $\sigma_{2k+1}(A)\neq 0$ for some $k\ge 1$, extending original Weyl's formula (for $\sigma_{2k}(A)$) in \cite{Weyl}.

Intrinsicality is related to rigidity of immersed hypersurfaces, and we refer to \cite{S, SP} and references therein for results of local and global rigidity of hypersurfaces in space forms. In view of Lemma \ref{l1}, hypersurface $\left(
M^{n},g\right)  \subset\left(  \mathbb{N}^{n+1}, \bar g\right)  $ is rigid in $ \mathbb{N}^{n+1}$ if $ M^{*}$ is connected, where 
\[M^{*}=\{x\in M \ |\  \sum_{k=0}^{[\frac{n}2]}\sigma^2_{2k+1}(A(x))>0, \ \mbox{rank}\  A(x)\ge 3\}.\]

From the proof of Lemma \ref{l3}, two embedded analytic hypersurfaces with the same induced metric are congruent. The global intrincality in Lemma \ref{l3} (with the outward orientation)  does not imply global rigidity in general, as indicated in the above figure, for example ( extendable to any dimension). 
\end{remark}

\begin{figure}[ptb]
\centering\includegraphics[
height=1.45in,
width=2.7in
]
{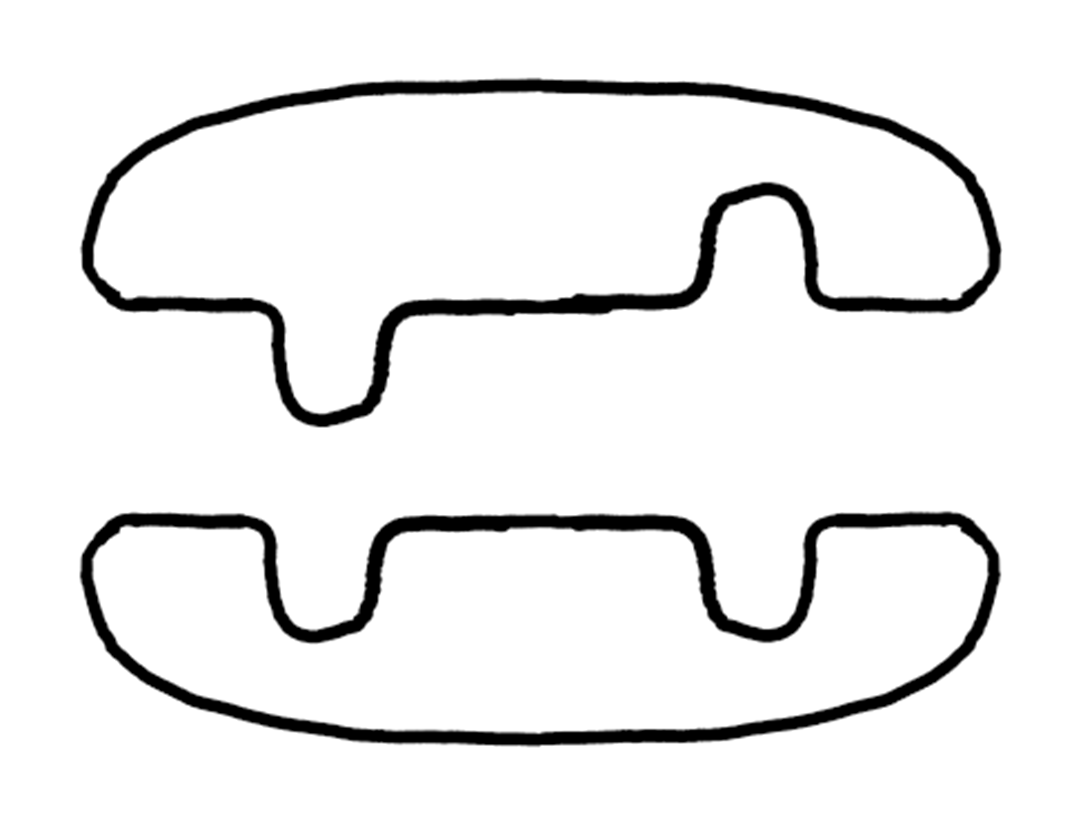}
\end{figure}

\textbf{Acknowledgments. }{The research of the first author is supported in part by an NSERC Discovery grant, the research of the second author is supported in part by an NSF grant.}

\end{document}